\documentclass[12pt]{article}
\usepackage{amsmath}
\usepackage{graphicx,psfrag,epsf}
\usepackage{enumerate}
\usepackage{natbib}
\usepackage{url} 
\usepackage{extarrows}
\usepackage{arydshln}

\newcommand{\blind}{0}

\newcommand{\bA}{\mbox{\bf A}}

\def\var{{\mbox{var}}}
\def\cov{{\mbox{cov}}}

\begin{document}

\def\spacingset#1{\renewcommand{\baselinestretch}%
{#1}\small\normalsize} \spacingset{1}


\if0\blind
{
  \title{\bf Better understanding of the multivariate hypergeometric distribution with implications in design-based survey sampling}
  \author{Xiaogang DUAN\thanks{
    The authors gratefully acknowledge \textit{the National Natural Science Foundation of China (No. 11771049)}}\hspace{.2cm}\\
    Department of Statistics, Beijing Normal University, Beijing  {\rm 100875}, China\\
    }
  \maketitle
} \fi

\if1\blind
{
  \bigskip
  \bigskip
  \bigskip
  \begin{center}
    {\LARGE\bf Better understanding of the multivariate hypergeometric distribution with implications in design-based survey sampling}
\end{center}
  \medskip
} \fi

\bigskip
\begin{abstract}

Multivariate hypergeometric distribution arises frequently in elementary statistics and probability courses, for simultaneously studying the occurence law of specified events, when sampling without replacement from a finite population with fixed number of classification. Covariance matrix of
this distribution is well known to be identical to its multinomial counterpart multiplied by $1-(n-1)/(N-1)$, with $N$ and $n$ being population and sample sizes, respectively. It appears to however, have been less discussed in the literature about the meaning of this relationship, especially regarding the specific form of the multiplier.
Based on an augmenting argument together with probabilistic symmetry, we present a more transparent understanding for the covariance structure of the multivariate hypergeometric distribution.
We discuss implications of these combined techniques and provide a unified description about the relative efficiency for estimating population mean  based on simple random sampling, probability proportional-to-size sampling and adaptive cluster sampling, with versus without replacement.
We also provide insight into the classic random group method for variance estimation.

\end{abstract}

\noindent%
\vfill

\newpage
\spacingset{1.45} 

\section{Introduction}\label{sec:intro}

%

Consider a finite population with $N$ units classified into $K$ subgroups each of population size  $N_k$ for $k=1,\ldots,K$.
Denote $a_k$ as the number observed falling into the $k$th category within $n$ sequential draws without replacement, one unit each draw, and $b_k$ as the corresponding number with replacement.
It is known that $(a_1,\ldots,a_K)^T$ follows a multivariate hypergeometric distribution, $(b_1,\ldots,b_K)^T$ follows a multinomial distribution,
and meanwhile for $k,l=1,\ldots,K$,
\begin{eqnarray}\label{intro-basic:relation}
\cov(a_k,a_l) &=& \cov(b_k,b_l) \left(1-\frac{n-1}{N-1}\right) \nonumber\\
&=&
\left\{
\begin{aligned}
n\frac{N_k}{N}\left(1-\frac{N_k}{N}\right)\left(1-\frac{n-1}{N-1}\right), &\quad k=l,\\
-n\frac{N_k}{N}\frac{N_l}{N}\left(1-\frac{n-1}{N-1}\right), &\quad k\neq l.
\end{aligned}
\right.
\end{eqnarray}

A basic solution towards equation (\ref{intro-basic:relation}) is to introduce an $n\times K$ matrix of indicator variables, each row recording an outcome for its associated draw, and express $a_k$ and $b_k$ as column sums of this indicator matrix. The desired result follows from an argument of covariance form of any single as well as any pair of indicator variables.

The above $n\times K$ matrix is observable in real sampling practice, but appears not enough for fine understanding results in equation (\ref{intro-basic:relation}).
In next section, we march ahead by augmenting this observable matrix, by an additional unobservable $(N-n)\times K$ of imagined indicator matrix which record possible realization if we would be willing to take $N$ instead of $n$ sequential draws without replacement. Based on this augmented matrix and its symmetric probability structure, we provide a more transparent depiction of components in equation (\ref{intro-basic:relation}), particularly for the multiplier $1-(n-1)/(N-1)$. It turns out that $(n-1)/(N-1)$ comes out because for any specific draw, $n-1$ additional draws are performed, each contributing a common negative covariance of order $-(N-1)^{-1}$.

Following this, we discuss implications of the combined techniques of augmenting and symmetrization in several context closely related to the fundamental simple random sampling in survey sampling textbook (e.g.: Cochran 1977; Thompson 2012). It turns out that the relative efficiency for estimating the population total, in the sense of without replacement relative to with replacement, are the common constant appeared in equation (\ref{intro-basic:relation}), for simple random sampling, probability proportional-to-size (PPS) sampling (Hansen \& Hurvitz, 1943) and the adaptive cluster sampling (Thompson 1990). We also provide an application of the above combined techniques in the classic random group method for variance estimation (Wolter 2007, chap 1). Finally, we give a brief conclusion.




\section{Augmenting, symmetrization, and their application}

For simple random sampling without replacement from a finite population of size $N$ with $K$ subgroups each of size $N_k$ for $k=1,\ldots,K$, we introduce the following $N\times K$ indicator matrix to finely record the entire sampling process
$$
\bA =
\begin{array}{c@{\hspace{-5pt}}l}
\left(\begin{array}{ccccc}
a_{1,1}& \cdots & a_{1,k} & \cdots & a_{1,K} \\
\vdots & \vdots & \vdots & \vdots & \vdots\\
a_{i,1}& \cdots & a_{i,k} & \cdots & a_{i,K} \\
\vdots & \vdots & \vdots & \vdots & \vdots\\
a_{n,1}& \cdots & a_{n,k} & \cdots & a_{n,K} \\ \hdashline[2pt/2pt]
a_{n+1,1}& \cdots & a_{n+1,k} & \cdots & a_{n+1,K} \\
\vdots & \vdots & \vdots & \vdots & \vdots\\
a_{N,1}& \cdots & a_{N,k} & \cdots & a_{N,K} \\ 
\end{array}\right)
&
\begin{array}{l}
\left.\rule{0mm}{17mm}\right\}{\tiny\mbox{Observable}}\\
\\
\left.\rule{0mm}{9.5mm}\right\}{\tiny\mbox{Unobservable}}
\end{array}
\end{array},
$$
where the binary variable $a_{i,k}$ takes the value $1$ if the $i$th random draw falls into the $k$th subgroup and zero otherwise, $i=1,\ldots,N; k=1,\ldots,K$.

Clearly, column sums of $\bA$ is the constant row vector $(N_1,\ldots,N_K)$. As a comparison, column sums of the first $n$ rows of $\bA$ is random and when transposed follows the multivariate hypergeometric distribution. Besides, two types of probabilistic symmetry hold readily from the definition of $\bA$:
\begin{itemize}
\item[P1.] All the row vectors of $\bA$ are commonly distributed.
\item[P2.] For any $N\times 2$ submatrix of $\bA$, all the $N(N-1)$ pair of entries selected from this submatrix, each pair consisting entries from a different column but not within the same row, are commonly distributed.
\end{itemize}





%

Let $a_k=\sum_{i=1}^n a_{i,k}$ for $k=1,\ldots,K$. The random vector $(a_1,\ldots,a_K)^T$ follows a multivariate hypergeometric distribution. We first show the variance of $a_k$ for fixed $k$, then derive the covariance of any pair of $a_k$ and $a_l$ for $k\neq l$.

By the row-wise probabilistic symmetry of $\bA$, or equivalently property P1, it follows
$$
\begin{aligned}
\var(a_k) &= \sum_{i=1}^n \var(a_{i,k}) + \sum_{i=1}^n \sum_{j\neq i} \cov(a_{i,k},a_{j,k})\\
&= n\{\var(a_{1,k}) + (n-1)\cov(a_{1,k},a_{2,k})\}.
\end{aligned}
$$
On the other hand, by the symmetry property P2, the following equations hold
$$
\begin{aligned}
\cov(a_{1,k},a_{2,k}) &= \cov(a_{1,k},a_{3,k}) = \cdots = \cov(a_{1,k},a_{N,k}) \\
&= \frac{1}{N-1}\cov\left(a_{1,k},\sum_{i=1}^N a_{i,k} - a_{1,k}\right)\\
&= \frac{1}{N-1}\cov\left(a_{1,k},N_k - a_{1,k}\right) = -\frac{1}{N-1}\var(a_{1,k}).
\end{aligned}
$$
Therefore, for sampling without replacement, we have
$$
\begin{aligned}
\var(a_k) = n\var(a_{1,k})\left(1-\frac{n-1}{N-1}\right).
\end{aligned}
$$

The above derivation indicates that, while both sampling with and without replacement share the common variance $n\var(a_{1,k})$, the additional covariance terms for sampling without replacement comes out because of $n-1$ extra draws each negatively correlated with the first draw with a common covariance equal to $-\var(a_{1,k})/(N-1)$.

For covariance of $a_{k}$ and $a_{l}$, following similar rules, we have
$$
\begin{aligned}
\cov(a_k,a_l) &= \sum_{i=1}^n \cov(a_{i,k},a_{i,l}) + \sum_{i=1}^n\sum_{j\neq i} \cov(a_{i,k},a_{j,l})\\
&= n\left\{\cov(a_{1,k},a_{1,l}) + (n-1)\cov(a_{1,k},a_{2,l})\right\}\\
&= n\left\{\cov(a_{1,k},a_{1,l}) + \frac{n-1}{N-1}\cov\left(a_{1,k},N_{l}-a_{1,l}\right)\right\}\\
&= n\cov(a_{1,k},a_{1,l})\left(1-\frac{n-1}{N-1}\right).
\end{aligned}
$$

When the underlying sampling is with replacement, all the covariance in double summation of both $\var(a_k)$ and $\cov(a_k,a_l)$ disappear. Therefore, compared to sampling with replacement, without replacement reduces the covariance matrix by a constant proportion $1-(n-1)/(N-1)$. It turns out that the same proportionate constant appears several times in classic design-based survey sampling, which are described in more detail in the next section. 

\section{Implications in design-based survey sampling}

For design-based survey sampling, the target population consists of a finite collection of units $\{1,\ldots,N\}$ say, with usually known population size $N$, and each unit is associated with a fixed value $Y_i$ for $i=1,\ldots,N$. A popular study aim is to estimate the population mean $\bar{Y}=N^{-1}\sum_{i=1}^N Y_i$ based on a probability sample $\{y_1,\ldots,y_n\}$ of size $n$ according to some prescribed sampling design.

Simple random sampling, which assigns equal probability for all possible samples, is fundamental in classic survey sampling context. Deep insights into a complex sampling design usually result if it could be viewed as some kind of simple random sampling. Celebrated examples include PPS sampling and adaptive cluster sampling. It will be shown that the combined techniques of augmenting and symmetrization provide insights into sampling designs closely related to simple random sampling. We also describe an application of these techniques in understanding a result in random group method for estimating the population variance in simple random sampling.


\subsection{Simple random sampling}

For simple random sampling, we estimate the population mean by the sample mean $\bar{y}=n^{-1}\sum_{i=1}^ny_i$.
The sample mean depends on population units actually selected within $n$ sequential draws without replacement $(i_1,\ldots,i_n)$ say.
The actually selected $n$ ordered units is a segment of an augmented vector $(i_1,\ldots,i_n, i_{n+1},\ldots,i_N)$, which itself is a possible permutation of $(1,\ldots,n,n+1,\ldots,N)$. Clearly, all entries $i_j (j=1,\ldots,N)$ are commonly distributed, and all $N(N-1)$ pairs $(i_{j_1},i_{j_2}) (1\leq i_{j_1}\neq i_{j_2}\leq N)$ are commonly distributed. Therefore, $\var(n\bar{y})=n\{\var(y_{1}) + (n-1)\cov(y_{1},y_{2})\}$.
On the other hand, $\var(y_{1})=N^{-1}\sum_{i=1}^N(Y_{i}-\bar{Y})^2$, and
$$
\begin{aligned}
\cov(y_{1},y_{2}) 
&= \frac{1}{N-1}\cov\left(y_{1},\sum_{i=1}^N Y_{i} - y_{1}\right) = -\frac{1}{N-1}\var(y_{1}).
\end{aligned}
$$
Consequently, for simple random sampling without replacement, the variance of $\bar{y}$ is
$$
\begin{aligned}
\var(\bar{y}) = \frac{1}{nN}\sum_{i=1}^N(Y_{i}-\bar{Y})^2\left(1-\frac{n-1}{N-1}\right).
\end{aligned}
$$
As a comparison, ignoring the contribution $-(n-1)/(N-1)$ due to without replacement, we obtain the variance of sample mean for simple random sampling with replacement.

%

\subsection{PPS sampling}

In the presence of unit-level auxiliary information, it is usually more efficient to incorporate them into the underlying sampling strategy.
A well-known example is the PPS sampling with replacement together with the Hansen-Hurvitz estimator $\hat{Y}_{HH}^{wr}=n^{-1}\sum_{j=1}^n (y_{i}/z_{i})$, for estimating the population total with an auxiliary size variable $M_i$ taking value of positive integers. Here, $Z_i = M_i/t_M$ with $t_M=M_1+\cdots+M_N$ for $i=1,\ldots,N$, is the selecting probability for the $i$th population unit in a single draw from the population with replacement; $y_i$ and $z_i$ denote the values of $Y$ and $Z$ for $i$th selection, $i=1,\ldots,n$.

The Hansen-Hurvitz estimator under PPS sampling may be interpreted as a simple sample average of size $n$, based on random sampling with replacement from a finite population of size $t_M$, with population variables taking the values of $Y_{1}/Z_1$ for unit $1$ through $M_1$, $Y_{2}/Z_2$ for units  $M_1+1$ through $M_1+M_2$, etc. The mean of this locally extended population is the population total $t_Y$ of the original population, that is $t_M^{-1}\{M_1(Y_{1}/Z_{1}) + \cdots + M_N(Y_{N}/Z_{N})\}=Y_1+\cdots+Y_N=t_Y$. The variance of the extended population is $t_M^{-1}\{M_1(Y_{1}/Z_{1}-t_Y)^2 + \cdots + M_N(Y_{N}/Z_{N}-t_Y)^2\}$. By results in the preceding subsection,
$\var(\hat{Y}_{HH}^{wr}) = n^{-1}\sum_{i=1}^N Z_i(Y_{i}/Z_{i}-t_Y)^2$. Moreover, we could improve this strategy by simple random sampling from the extended population without replacement, and obtain a more efficient estimator which reduces $\var(\hat{Y}_{HH}^{wr})$ by a factor $1-(n-1)/(t_M-1)$.
%
%

\subsection{Adaptive cluster sampling}

Adaptive sampling is often encountered in spatial context (Thompson 1990). For adaptive cluster sampling, population units enter into the sample either through an initial or a subsequent sampling. In the absence of auxiliary information, initial samples are usually obtained by simple random sampling, with or without replacement. Neighbouring population units of each initially selected unit then fall into the sample, only if the initially selected unit satisfies certain prescribed conditions $C$; for example $C=\{y>y_0\}$ and $y_0$ is a fixed number. If any of these neighbouring units satisfies the condition, still more units will be added into the sampling process.

An essential concept for adaptive sampling is network, which is a subset of population units such that any one unit initially selected in this subset  leads to the entire collection entering into the final sample. All networks form a partition of the original population. In other words, the population $U=\{1,\ldots,N\}$ could be reorganized into $K (K\leq N)$ distinct networks $U^{(1)},\ldots,U^{(K)}$ with $U^{(k)}\cap U^{(l)}=\emptyset$ for any $k\neq l$. Compared to conventional sampling, one challenge is to construct estimators of population mean by making effective use of observed information for all sampled units. Thompson (1990) described an innovative procedure, which could be understood as a simple sample average based on simple random sampling, with or without replacement.

In fact, the original population could be locally flattened out int the sense that we associate each unit in a network, instead of its own $Y$ value,  a common value, the average of $Y$ values in this network. For this locally flattened population, its population mean equals clearly to the original population mean. The population variance is $N^{-1}\sum_{k=1}^K N_k(\bar{Y}^{(k)}-\bar{Y})^2$, where $\bar{Y}^{(k)}=N_k^{-1}\sum_{i\in U^{(k)}}Y_i$ is population mean for the $k$th network, and $N_k$ is its total number of population units. With these preparation, one can happily estimate the population mean by sample mean of initially selected units, not for their $Y$ values, but for their locally shared average of $Y$ values across the respective network in which they reside. Compared to initial sampling with replacement, the initial sampling without replacement possesses an relative efficiency gain of $1-(n_1-1)/(N-1)$, where $n_1$ is number of units initially selected.

%

\subsection{Random group method}


In the case of estimating the population mean with observations $\{y_1,\ldots,y_n\}$ obtained by simple random sampling from a finite population without replacement, the random group method splits the full sample into $K$ non-overlapping subgroups, each of fixed sample size $n_k$ for $k=1,\ldots,K$, and attempts to evaluate the variance of a full sample estimator based on the sample variance of the $K$ subgroup estimators.

Let $\overline{y}^{(k)}$ be the sample mean of observations in the $k$th subgroup for $k=1,\ldots,K$. It is known that the sample variance of $\overline{y}^{(1)},\ldots,\overline{y}^{(K)}$ is the simple average of all half squared pair differences $(\overline{y}^{(k)}-\overline{y}^{(l)})^2/2$ for $k\neq l$ and $k,l=1,\ldots,K$.  We need to evaluate the expectation of $(\overline{y}^{(k)}-\overline{y}^{(l)})^2$,
with the randomness from both the first stage sampling without replacement as well as the successive random grouping.

Note that
\begin{equation*}
\frac{1}{\binom{N}{n}\times\{\binom{n}{n_1}\binom{n-n_1}{n_2}\cdots 1\}} = \left\{\frac{N!}{n_1!n_2!\cdots n_K! (N-n)!}\right\}^{-1}.
\end{equation*}
This implies that probabilistic structure of the two-stage random grouping coincides with those generated from $n$ sequential sampling without replacement from the same population, and ascribing the first $n_1$ draws into the first subgroup, the subsequent $n_2$ draws into the second subgroup, and so on.

In other words, if we denote $(i_1,\ldots,i_n)$ as the $n$ sequentially selected result, the first subgroup consists of units $(i_1,\ldots,i_{n_1})=\omega^{(1)}$, and the second subgroup consists of units $(i_{n_1+1},\ldots,i_{n_1+n_2})=\omega^{(2)}$, etc. Clearly, all single units in the augmented vector $(i_1,\ldots,i_N)$ are symmetric in probability, and all pair units of $(i_1,\ldots,i_N)$ are commonly distributed. Consequently
$$
\begin{aligned}
E\left\{(\overline{y}^{(k)}-\overline{y}^{(l)})^2\right\} &= \var(\overline{y}^{(k)}) + \var(\overline{y}^{(l)}) - 2\cov(\overline{y}^{(k)},\overline{y}^{(l)})\\
&= \frac{1}{N-1}\sum_{i=1}^N(Y_{i}-\bar{Y})^2\left(\frac{N-n_k}{Nn_k}+\frac{N-n_l}{Nn_l}\right) - 2\cov(y_1,y_2) \\
&= \frac{1}{N-1}\sum_{i=1}^N(Y_{i}-\bar{Y})^2\left(\frac{N-n_k}{Nn_k}+\frac{N-n_l}{Nn_l}\right) + \frac{2}{N-1}\var(y_1) \\
&= \frac{1}{N-1}\sum_{i=1}^N(Y_{i}-\bar{Y})^2(n_k^{-1}+n_l^{-1}).
\end{aligned}
$$
This suggests that we can estimate the population variance unbiasedly by adjusting the sample variance of $\overline{y}^{(k)}, \ldots,\overline{y}^{(K)}$.

\section{Conclusion}

We have provided a more intuitive depiction for the covariance structure of the multivariate hypergeometric distribution, and explored its fine relationship with the covariance structure of the associated multinomial distribution. The basic techniques consist of augmenting as well as some argument about the probabilistic symmetry of the augmented sampling. Both techniques are easy to understand, given some knowledge of basic probability space in elementary probability courses.
The provided materials should be a valuable supplement to understanding efficiency comparisons for some important sampling designs.

\section*{References}

\begin{description}

%
%

\item

\item Cochran W.G. (1977). Sampling Techniques, 3rd ed. John Wiley \& Sons.

\item Hansen M.H. \& Hurwitz W.N. (1943). On the theory of sampling from finite populations. \textit{Ann. Math. Statist}, \textbf{14},333--362.

\item Thompson S.K. (1990). Adaptive cluster sampling. \textit{Jour. Amer. Statist. Assoc.}, \textbf{85}, 1050--1059.

\item Thompson S.K. (2012). Sampling, 3rd ed. John Wiley \& Sons.

\item Wolter KM (2007). Introduction to variance estimation. Springer.

\end{description}

\end{document}